\documentclass[reqno,12pt,english]{amsart}
\usepackage{amssymb,amsmath,amstext,amsthm,amsfonts,amscd}
\usepackage[dvips]{graphicx}
\usepackage{psfrag,euscript,a4wide,stmaryrd,wasysym}

\newtheorem{theorem}{Theorem}

\newtheorem{maintheorem}{Theorem}

\newcommand{\cmt}{\begin{maintheorem}}
\newcommand{\fmt}{\end{maintheorem}}

\newtheorem{maincorollary}[maintheorem]{Corollary}

\newcommand{\cmc}{\begin{maincorollary}}
\newcommand{\fmc}{\end{maincorollary}}

\newtheorem{T}{Theorem}[section]
\newcommand{\cte}{\begin{T}}
\newcommand{\fte}{\end{T}}

\newtheorem{Corollary}[T]{Corollary}
\newcommand{\cco}{\begin{Corollary}}
\newcommand{\fco}{\end{Corollary}}

\newtheorem{Proposition}[T]{Proposition}
\newcommand{\cpr}{\begin{Proposition}}
\newcommand{\fpr}{\end{Proposition}}

\newtheorem{Lemma}[T]{Lemma}
\newcommand{\cle}{\begin{Lemma}}
\newcommand{\fle}{\end{Lemma}}

\newcommand{\csle}{\begin{Lemma}}
\newcommand{\fsle}{\end{Lemma}}

\newtheorem{Remark}[T]{Remark}
\newcommand{\cre}{\begin{Remark}}
\newcommand{\fre}{\end{Remark}}

\newtheorem{Definition}[T]{Definition}
\newcommand{\cde}{\begin{Definition}}
\newcommand{\fde}{\end{Definition}}

\newcommand{\SF}{{\mathcal F}}

\newcommand{\SK}{{\mathcal K}}

\newcommand{\SM}{{\mathcal M}}

\newcommand{\SP}{{\mathcal{P}}}

\newcommand{\SR}{{\mathcal R}}

\newcommand{\ST}{{\mathcal{ T}}}
\newcommand{\SU}{{\mathcal U}}

\newcommand{\dem}{\begin{proof}}
\newcommand{\cqd}{\end{proof}}

\newtheorem{proposition}[theorem]{Proposition}
\newtheorem{lemma}[theorem]{Lemma}

\theoremstyle{definition}
\newtheorem{remark}[theorem]{Remark}

\newtheorem{definition}[theorem]{Definition}

\newcommand{\pf} {{\bf Proof: }}
\newcommand{\field}[1]{\mathbb{#1}}

\newcommand{\ov} {\overline}

\renewcommand{\natural}{\field{N}}

    \newcommand{\Ga}{\Gamma}

      \newcommand{\La}{\Lambda}

\newcommand{\vsi}{\varsigma}

\newcommand{\vr}{\varphi}

\newcommand{\NE}{\operatorname{{\text{\em NEnd}^1}}}
\newcommand{\DI}{\operatorname{{\text{\em Diff}^1}}}

\newcommand{\diam}{\operatorname{{diam}}}

\begin{document}


\author{Armando Castro}
\address{Departamento de Matem\'atica, Universidade Federal da Bahia\\
Av. Ademar de Barros s/n, 40170-110 Salvador, Brazil.}
\email{armandomat@pesquisador.cnpq.br}

\subjclass[2000]{Primary 37C50, 37C20; Secondary 37D05, 37C75}

\date{\today}

\keywords{Ergodic Theory, Structural Stability Conjecture for
Endomorphisms, Ergodic Closing Lemma, Nonsingular Endomorphisms}

\thanks{Work carried out at the  Federal University of
Bahia. Partially supported by CNPq (PQ 10/2007) and UFBA}

\title[Ergodic Closing Lemma]{The Ergodic Closing Lemma for nonsingular Endomorphisms}

\maketitle

\begin{abstract}
We generalize Ma\~n\'e's Ergodic Closing Lemma to the context of
$C^1$-Endomorphisms without singularities.
\end{abstract}

\section{Introduction}

Let $M$ be a finite dimensional compact Riemannian manifold and let
$\NE(M)$   denote the set of $C^1$ {\em nonsingular endomorphisms}
defined in $M$. By $g:M \to M$ to be a {\em nonsingular
endomorphism} we mean that the derivative $Dg(p)$ of $g$ in each
point $p \in M$ is a linear isomorphism. We endow $\NE(M)$ with the
$C^1$ topology, and denote its corresponding metrics by $d_1$;
therefore $\NE(M)$ is an open subset of the complete space $C^1(M)$
whose elements are $C^1-$endomorphisms defined in $M$.

The main purpose of this paper is to prove

\cmt \label{theo1} {(Ergodic Closing Lemma for nonsingular
endomorphisms.)} Let $M$ be a compact manifold. Then,there exists a
residual subset $\SR \subset \NE(M)$ such that for any $f \in \SR$,
the set of $f-$invariant probabilities $\SM_1(f)$ is the closed
convex hull of ergodic measures supported on periodic orbits of $f$.
\fmt

For the proof of theorem above, we also prove another version of
Ergodic Closing Lemma for endomorphisms (see Th. \ref{theop}) in the
next section.

\begin{remark}
The Ergodic Closing Lemma is a result about shadowing by periodic
orbits. Although the classical notion of shadowing is not generic
even among $C^1$-diffeomorphisms (see \cite{DA}), Ergodic Closing
Lemma asserts that $C^1-$generically most orbits in a
measure-theoretical point of view can be shadowed by periodic
orbits. Using ideas in \cite{YC}, \cite{COP}, \cite{C2}, and the
Ergodic Closing Lemma we also obtained a new criteria of generic
Hyperbolicity/Expansion based in periodic sets \cite{C3}.
\end{remark}

{\bf Acknowledgements: }I am grateful to professors Marcelo Viana,
Vilton Pinheiro, and Paulo Varandas for conversations on topics
related to this paper. I am also grateful to my beloved wife Maria
Teresa Gilly, whose shadow shadows mine.

\section{Proof of the Ergodic Closing Lemma for nonsingular endomorphisms}

\label{sec3} Let us start by fixing some notation. Given $x \in M$,
we define $B_\epsilon(f, x)$ as an $\epsilon-$neighborhood of the
orbit of $x$. Define  $\Sigma(f)$ as the set of points $x \in M$
such that for every neighborhood $\SU$ of $f$ and every $\epsilon>
0$, there exist $g \in \SU$ and $y \in M$ such that $y \in Per(g)$,
$g= f$ on $M \setminus B_\epsilon(f, x)$ and $d(f^j(x), g^j(y)) \leq
\epsilon$, $\forall 0 \leq j \leq m$, where $m$ is the $g-$period of
$y$.

The {\bf residual
} version of the Ergodic Closing Lemma (our Th. \ref{theo1}) is a consequence of
the following result:
\begin{theorem}\label{theop}
For any nonsingular endomorphism $f$, $\Sigma(f)$ is a total probability set,
that is, $\Sigma(f)$ is a full probability set for any $f-$invariant
probability.
\end{theorem}

\begin{remark}
Theorem \ref{theop} above, with $f \in \DI(M)$ instead of $f \in
\NE(M)$ in its statement, was the former Ergodic Closing Lemma
proved by Ma\~n\'e. In fact, Ma\~n\'e, in \cite{Ma82}, did not
explicitly give the proof of his corresponding classical residual
version, although he stated such version in \cite{Ma87}. Recently,
while we were writing this article,  Abdenur et al in \cite{ABC}
filled out this gap.
\end{remark}

\begin{definition}{($\epsilon-$shadowing by a periodic point.)} Let $f$ and $g$ maps on
a compact metric space $\La$. Given $\epsilon> 0$ and $x \in \La$, we say that
a $g-$periodic point $p$ with period $n$  {\em $\epsilon-$shadows $x$} iff
$d(g^j(p), f^j(x))< \epsilon, \forall 0 \leq j \leq n$.
\end{definition}

\begin{definition} Let $\epsilon> 0$ and a neighborhood $\SU \ni f$, $\SU \subset \NE(M)$.
 We define {\em $\Sigma(f, \SU, \epsilon)$} as
the set of points $x \in M$ such that there exist $g \in \SU$ and $y \in
M$ such that $y \in Per(g)$, $g= f$ on $M \setminus B_\epsilon(f,
x)$ and $d(f^j(x), g^j(y)) \leq \epsilon$, $\forall 0 \leq j \leq
m$, where $m$ is the $g-$period of $y$. That is, $\Sigma(f, \SU, \epsilon)$ is
the set of points $x \in M$ which are $\epsilon-$shadowed by a
periodic point $y \in Per(g)$, for some $g \in \SU$.
Everytime there is no chance of misunderstanding, we will just
write $\Sigma(\SU, \epsilon)$ instead of $\Sigma(f, \SU, \epsilon)$.
\end{definition}

If we take a nested neighborhood basis $\SU_n$ of $f$ in $\NE(M)$ then
$$
\Sigma(f)= \cap_{n \in \natural} \Sigma(f, \SU_n, 1/n).
$$

Therefore, Th. \ref{theop} is an immediate consequence of

\begin{proposition}\label{prop}
For any nonsingular endomorphism $f$, any neighborhood $\SU$ of $f$
and $\epsilon> 0$, $\Sigma(\SU, \epsilon)= \Sigma(f, \SU,\epsilon)$
is a total probability set for $f$.
\end{proposition}

The proof of Proposition \ref{prop} is quite long, and it is a
consequence of lemmas and Theorems we prove in the sequel. Such
proof has  two main parts.

Part 1 consists in an improvement of Closing Lemma (see \cite{P1},
\cite{P2}, \cite{PR}, \cite{LW1}), stating that given a nonsingular
endomorphism
 $f$, for any point $x \in M$ that returns sufficiently close to itself,
there is an iterate $y= f^{m(x)}(x)$ such that we can perturb $f$
into a $g$ for which there exists a periodic orbit that shadows $y$.
The precise statement of this part corresponds to lemma
\ref{lesomb}, whose proof we write down further in this paper. Note
that this part is entirely topological, and it does not use
measure/ergodic theoretical arguments.

Part 2 uses Birkhoff Theorem (ergodicity) and a Vitali's Covering
argument to prove that the set of points that are shadowed by a
periodic point of some nearby endomorphism $g$ has total probability
for $f$. The core of this part consists specially in Lemma
\ref{lerg} and Th. \ref{teou}.

For the statement of the Perturbation lemmas that we will need for the
 proof of Lemma \ref{lesomb}, let us introduce some
notation. As the manifold $M$ is compact, there is $\delta$ such
that $\{\exp_p, p \in M\}$ is an equilipschitz family of
diffeomorphisms, such that each exponential map $\exp_p$ embeds
$B(0, \delta)$ in a neighborhood $B_p$ of $p$. Given $p \in M$, we
define a metric $d'= d'_p:B_p \times B_p \to [0, +\infty)$ given by
$$
d'(x, y):= |\exp^{-1}_p(x)- \exp^{-1}_p(y)|.
$$

 Obviously, $d'$ is Lipschitz-equivalent to the manifold usual
metric restricted to $B_p$. Setting $d'$ as the metric in $B_p$,
then $\exp_p$ isometrically maps $B(0, \delta)$ on $B_p= B'(p,
\delta)$, where the quote ' signs the ball in the metric $d'$.

\begin{lemma}{\cite{PR}}\label{lelift}
For any $\eta > 0$, there is an $\alpha > 0$ such that for any $f
\in \NE(M)$, any $q \in M$, any two points $v_1, v_2 \in T_qM$ with
$B(v_2, |v_1- v_2|/\alpha) \subset B(0, \delta) \subset T_qM$, there
is a diffeomorphism $h= h_{q, \alpha, v_1, v_2}: M \to M$, called an
$\alpha-$lift, such that:
\begin{enumerate}
\item $h(\exp_q(v_2))= \exp_q(v_1)$;
\item The closure of set of points where $h$ differs to the identity
is contained in $\exp_q(B(v_2, |v_1- v_2|/\alpha)$;
\item $d_1(h f, f)< \eta$.
\end{enumerate}
\end{lemma}

\begin{definition}{(Dynamical neighborhood.)}
We say that a neighborhood $V$ of a point $p \in M$ is  $N-${
dynamical} for $f$ if each connected component $\cup_{j= 0}^N
f^{-j}(V)$ contains exactly one point of $\cup_{j= 0}^N
f^{-j}(\{p\})$.
\end{definition}

\begin{lemma}{\cite{F}, \cite{LW1}}\label{lelinea}
Let $f \in \NE(M)$, $p \in M$, $N \geq 1$ given such that all terms in
$\cup_{j= 0}^{N+ 1} f^{-j}(p)$ are distinct. Then, for any $\eta>
0$, there is a $\beta> 0$, and a map $f_1 \in \NE(M)$, called a local
linearization of $f$ with the following properties (1)-(5).
\begin{enumerate}
\item
$B'(p, \beta)$ is  $(N+ 1)$-{\em dynamical } for both $f$ and $f_1$,
and $f^{-j}(B'(p, \beta))= f_1^{-j}(B'(p, \beta))$ for $j= 1, \dots,
N+1$.

\item For $q \in \cup_{j= 1}^{N+ 1} f^{-j}(p)$, let $V(q)$
be the open connected component of $\cup_{j= 1}^{N+ 1} f^{-j}(B'(p,
\beta/4))$ containing $q$. Then, $f_1|_{V(q)}= \exp_{f(q)} \circ
(T_q f) \circ \exp_q^{-1}$.

\item $f_1^{N+1}(x)= f(x)$, $\forall x \in f^{-N-1}(B'(p,
\beta))$.

\item $f_1= f$ on $M \setminus \cup_{j= 1}^{N+ 1} f^{-j}(B'(p,
\beta))$.

\item $d_1(f_1, f)< \eta$.

\end{enumerate}

\end{lemma}

\begin{remark} \label{rent}
For the sequel, we need to emphasize two trivial, but important
consequences of the last lemma. The first one is that if for some $k
\in \natural$ we have that $x, f^k(x)$ are both out of $\cup_{j=
0}^{N} f^{-j}(B'(p, \beta))$, then $f^k(x)= f_1^k(x)$. The second
one is that, as $\eta$ goes to $0$, one can take $\beta$ arbitrarily small in
the last lemma.
\end{remark}

\begin{theorem}{(Th. A in  \cite{LW1}.)} \label{teoalw}
Let $(\ST, T_q)$ a complete tree of isomorphisms associated to the
pre-orbit of a point $q_0 \in M$, that is, a collection of
$n-$dimensional inner product spaces $E_{q}$ and isomorphisms $T_q:
E_q \to E_{q_0}$ associated to each $q$ in the pre-orbit of $q_0$,
with $T_{q_0}$ equal to identity. Given $\alpha> 0$, there are
$\rho> 0$ and $N \geq 1$ such that: for any ordered set $X= \{x_0
\prec \dots \prec x_t \} \subset E_{q_0}$, there is a point $y \in X
\cap B(x_t, \rho |x_0- x_t|)$ such that for any branch $\Ga= \{q_0,
q_1, \dots, \}$ of $\ST$, there is a point $w= w(\Ga) \in X \cap
B(x_t, \rho |x_0- x_t|)$ which is before $y$ in the order of $X$,
together with $N+ 1$ points $c_0(\Ga)= c_0, \dots, c_N(\Ga)= c_N \in
B(x_t, \rho |x_0- x_t|)$ satisfying the following two conditions:
\begin{itemize}
\item $c_0= w, c_N= y$; and

\item $|T_{q_n}^{-1}(c_j)- T_{q_n}^{-1}(c_{j+1})| \leq \alpha
d(T_{q_n}^{-1}(c_{j+1})- T_{q_n}^{-1}(A))$, where $A:= \{x \in X,
w\prec x \prec y\} \cup \partial B(x_t, \rho |x_0- x_t|)$.

\end{itemize}
\end{theorem}

The next lemma is the main target in the first part of our
Ergodic Closing Lemma. It is a topological result, and it has nothing to do with
 Measure/Ergodic theory. It implies in particular that, given $\epsilon> 0$ and
  {\em any} $f-$recurrent
point $x$, then $x$ has an iterate which is $\epsilon-$shadowed by a periodic point of some
$g$ close to $f$.

\begin{lemma}
\label{lesomb} Given $f \in \NE(M)$, $p \in M$, $\epsilon> 0$ and a
neighborhood $\SU$ of $f$, there exist $r> 0$, $\rho'> 1$ such that
if for some natural $t> 0$, we have $x, f^t(x) \in B'_{\ov r}(p)$,
with $0< \ov r \leq r$, then there exist $0 \leq t_1 < t_2 \leq t$
and $g \in \SU$ such that:
\begin{itemize}
 \item $w= f^{t_1}(x), y= f^{t_2}(x) \in \ov{B'_{\rho' \ov r}(p)}$;
 \item $g^{t_2- t_1}(w)= w$;
 \item $g(z)= f(z)$ for $z \notin B_\epsilon(f, x)$ and $d(g^j(w), f^j(w)) \leq
\epsilon$, $\forall 0 \leq j \leq t_2- t_1$.
\end{itemize}
\end{lemma}
\pf Take an $\eta> 0$ such that the $\eta-$ball with center $f$ is
contained in $\SU$.
 Take $1> \alpha> 0$ such that $d_1(h \circ f, f) < \eta/2$, for any
$\alpha-$lift $h$.

Without loss of generality, we assume that $\epsilon< \alpha^2$. We
also assume that $\epsilon< \delta$.

We assume that $p$ is not periodic for $f$, otherwise, there is
nothing to prove.

This implies that all points in the pre-orbit of $p$ are distinct.
Let $\rho> 2$ and $N \geq 1$ be the numbers provided  by th.
\ref{teoalw}, for $\alpha> 0$ taken as above, and for $q_0= p$, each
$q_j$ to be some $j-$pre-image of $p$,  $E_{q_j}= T_{q_j} M$ and
$T_{q_j}= Df^j(q_j)$.

So, take $r> 0$ such that $r< \epsilon/(6\rho)$ and
$\diam(f^{-j}(B'(p, 3\rho r))< \epsilon$, $\forall j= 0, \dots,
N+1$. We assume that each connected component of $\cup_{j= 0}^{N+ 1}
f^{-j}(B'(p, 3\rho r))$ contains exactly one point $q_j \in
f^{-j}(p), j= 1, \dots, N+1$. In particular, if $z, f^{\hat t}(z) \in
B'(p, 3\rho r)$, then $\hat t> N+ 1$.

Now, assuming that $x, f^t(x) \in B'(p, \ov r)$, for some $0< \ov r
< r$ we can apply th. \ref{teoalw} to the set $X= \{x, f(x),\dots,
f^t(x)\} \cap B'(p, 3\rho \ov r)$ endowed with the order given by
the iterate number: if $f^k(x), f^{\hat k}(x) \in X$, then $f^k(x)
\prec f^{\hat k}(x) \Leftrightarrow k< \hat k$. Therefore, set
$\rho'= 3 \rho$. We then obtain $f^{t_2}(x)= y \in \{x, \dots,
f^t(x)\} \cap B'(f^t(x), \rho \cdot d'(f^t(x), x)) \subset B'(p,
\rho' \ov r)$ such that for any branch $\Ga= \{p= p_0, p_1, \dots,
p_n, \dots\}$ of the pre-orbit of $p$, there is $w= w(\Ga)=
f^{t_1}(x) \in \{x, \dots, f^t(x)\} \cap B'(f^t(x), \rho \cdot
d'(f^t(x), x))$, with $t_1= t_1(\Ga)< t_2$ together with points
$c_0= c_0(\Ga), \dots, c_N= c_N(\Ga) \in B'(f^t(x), \rho \cdot
d'(f^t(x), x))$ such that:

\vspace{0.1cm}

{(a)} $c_0= w, c_N= y$; and

{(b)} $|T_{p_j}^{-1}(c_j)- T_{p_j}^{-1}(c_{j+1})| \leq \alpha
d(T_{p_j}^{-1}(c_{j+1}), T_{p_j}^{-1}(A))$, where $A:= \{f^j(x) \in
X; t_1< j< t_2 \} \cup \partial B'(f^t(x), \rho \cdot d'(f^t(x),
x))$.

\vspace{0.2cm}

As $w= w(\Ga)$ and $y$ are both in $X$, there is a natural number
$k(\Ga) \geq 1$ such that $f^{k(\Ga)}(w(\Ga))= y$. Note that
$k(\Ga)> N+ 1$, as $\cup_{j= 0}^{N+1} f^{-j}(\{y\}) \cap B'(p, 3\rho
r)= y$, from our choice of $r$. Setting $z:= f^{k(\Ga)- N-
1}(w(\Ga))$, we see that $z$ does not depend on the branch $\Ga$ of
$p$, since $w(\Ga)$ and $y$ are in $X$, $f^{k(\Ga)}(w(\Ga))= y$ and
$y, N$ do not depend on $\Ga$. By our choice of $r$, since $y \in
B'(p, 3\rho r)$ there is a unique connected component $V_{N+1}
\subset f^{-(N+1)}(B'(p, 3\rho r))$ such that $z \in V_{N+1}$. Also,
there is a unique $p_{N+1} \in f^{-(N+1)}(p) \cap V_{N+1}$. From now
on, we fix $\Ga$ as some branch of $p$ containing $p_{N+1}$ (That
is, $\Ga= (p= p_0, \dots, p_{N+1}, \dots)$), and we consider all
constants $w, c_0, \dots, c_N$, $k$ obtained by applying th.
\ref{teoalw} with respect to such branch. For each $p_j, j= 0,
\dots, N-1$, let $h_{p_j}$ be the $\alpha-$kernel lift obtained by
treating in lemma \ref{lelift} $q= p_j$, $v_1=
[Df^j(p_j)]^{-1}(c_j)$, $v_2= [Df^j(p_j)]^{-1}(c_{j+1})$. Defining a
map $g: M \to M$ by
$$
g:= \left\{
\begin{matrix}
h_{p_j} \circ f_1 \text{ on } V(p_{j+1}); \\
f_1    \text{ on the rest of } M,
\end{matrix}
\right.
$$
we have that $g \in \NE(M)$ and $d_1(g, f)< \eta$. Thus $g \in \SU$.

Due to condition (b) above, the $g-$orbit from $w$ to $z$ never
touches the region in which $g \neq f_1$. Therefore, $g^{k-
(N+1)}(w)= f_1^{k- (N+1)}(w)$. By remark \ref{rent}, we also have
that $f^{k- (N+1)}(w)= f_1^{k- (N+1)}(w)$, and thus
$$
g^{k- (N+1)}(w)= f^{k- (N+1)}(w)= z.
$$
Now, it is easy to see that $g^{N+1}(z)= w$ and then $g^k(w)= w$. In
fact, $f_1^{N+1}(z)= y$, and the lifts $h_{p_{N-1}}, \dots h_0$
gradually and slightly modifies $f_1-$orbit segment joining $z$ and
$y$, in such way that $g^{N+1}(z)= w$ and $d(g^j(z), f^j(z)) \leq
d(g^j(z), f_1^j(z))+ d(f_1^j(z), f^j(z))< \epsilon$, $\forall j= 1,
\dots, N+1$.

\qed

Now we begin the second part of the proof of Proposition \ref{prop}.
Although   the main idea of this part is borrowed from \cite{Ma82},
the proofs we have written are presented in an abstract setting for
future use and bookkeeping purposes. This will also clarify the sort
of arguments which are used. Such arguments are basically measure
theoretical, and ergodic tools.

We start by introducing some notation. We say that a subset $C$ of the torus $T^s$ is a cube
if it can be written as $A= I_1 \times \dots \times I_s$, where the sets $I_i$ are
intervals of same length in $S^1$ (containing both, none, or one of its boundary points).
If $p_i$ is the middle point of $I_i$, we say that the point $(p_1, \dots, p_s)$
is the center of $A$. The length of the intervals $I_i$ is called the side of the cube.
For each $k \in \natural^+$, let $(\SP^{(k)}_j)_{j \in \natural^+}$ be a  sequence of partitions
of $T^s$ by cubes whose side is $2\pi/k^j$. For every atom $P$ of a partition $\SP^{(k)}_j$, we can associate
cubes $\hat P$ and $\tilde P$ having the same center of $P$, but with sides $2\pi/k^{j-1}$
and $6\pi/k^{j-1}$, respectively. If $x \in T^s$, denote by $P^{(k)}_j(x)$
the atom of $\SP^{(k)}_j$ containing $x$. Suppose that $M$ is isometrically embedded in $T^s$.

We recall the following useful fact on such kind of partitions:

\begin{lemma}\label{leprp}
For every probability measure $\mu$ on the Borel sets of $T^s$, every $\delta> 0$ and for all
odd natural $k$, the following inequalities holds for any $j \geq 1$:
$$
\mu(\{x; \mu(P^{(k)}_j(x)) \geq \delta \mu(\hat P^{(k)}_j(x))\}) \geq 1- \delta k^s
$$
and
$$
\mu(\{x; \mu(P^{(k)}_j(x)) \geq \delta \mu(\tilde P^{(k)}_j(x))\}) \geq 1- \delta 3^sk^s.
$$
\end{lemma}
\pf Done in \cite{Ma82}.

\qed

Let $f \in \NE(M)$, $\epsilon> 0$, a neighborhood $\SU$ of $f$
and an ergodic $\mu \in \SM_1(f)$ be given. Extend $\mu$ to a measure on $T^s$ by $\mu(A):= \mu(A \cap M)$,
for all Borel set $A \subset T^s$.

Let $\mho \subset M$ be some Borelian set and suppose that  $\mho(
r, \rho)$, where $r> 0$, $\rho> 1$, is some Borelian set whose
elements are points $x \in M$ with the following property:  if $y
\in B'_{r'}(x)$ for some $0< r' \leq r$ and $f^t(y) \in B'_{r'}(x)$,
for some $t> 0$ then there exist $0 \leq t_1< t$, such that
$f^{t_1}(y) \in \ov{B'_{\rho r'}(x)} \cap \mho$. { Take $r_i> 0$,
$\rho_l > 1$ two monotone sequences} converging respectively to 0
and $+\infty$.

Our first target in this second part  is to obtain an abstract
result (Th. \ref{teou})  which will be essential in both  proofs of
Proposition \ref{prop} and Th. \ref{theo1}.
Such result says that, if $\cup_{i, l} \mho(r_i, \rho_l)= M$, then
$\mho$ has total probability for $f$.

\begin{remark}
All results from this point of the paper up to Th. \ref{teou} do not use much regularity of $f$. In fact, specifically for the statements from Lemma \ref{lerg} up to Th. \ref{teou}, we only request  $f:M \to M$ to be
a Borelian map such that $\SM_1(f) \neq \emptyset$. This occurs, for instance, if $f$
is a continuous map.
\end{remark}

For each pair $(i, l)$, we can find
and odd natural $k= k(i, l)$ and $j(i, l)$ such that $\forall j \geq j(i, l)$ and $x \in T^s$ there exists $0 \leq r \leq r_i$ satisfying
$$
P^{(k)}_j(x) \subset B_r(x)
$$
and
$$
\hat P^{(k)}_j(x) \supset \ov{B_{\rho_l r}(x)},
$$
the balls here are taken in the torus.

The next lemma is where the $\mu$-ergodicity is necessary for the proof of Proposition \ref{prop}:

\begin{lemma}\label{lerg}
If \  $x \in  \mho( r_i, \rho_l)$, $j \geq j(i, l)$, $k= k(i, l)$ and $\mu(P^{(k)}_j(x)) \geq \delta \mu(\hat P^{(k)}_j(x))$, we have:
$$
\mu(\hat P^{(k)}_j(x) \cap \mho) \geq \delta \mu(\hat P^{(k)}_j(x)).
$$
\end{lemma}
\pf
As $\mu$ is ergodic, for $\mu-$typical $y \in M$, we have that
$$
\mu(\hat P^{(k)}_j(x) \cap \mho)= \lim_{n \to +\infty} \frac 1 n \#\{1 \leq t \leq n; f^t(y) \in \hat P^{(k)}_j(x) \cap \mho \},
$$
and
$$
\mu(P^{(k)}_j(x))= \lim_{n \to +\infty} \frac 1 n \#\{1 \leq t \leq n; f^t(y) \in P^{(k)}_j(x) \}.
$$
By the definition of $\mho(r_i, \rho_l)$, between any pair of natural numbers $n_1$ and $n_2$ such that $f^{n_1}(y), f^{n_2}(y) \in P^{(k)}_j(x) \subset B'_r(x)$, there exists $n_1\leq t_1< n_2$, such that $f^{t_1}(y) \in (\ov{B'_{\rho_l r}(x)} \cap \mho) \subset (\hat P^{(k)}_j(x) \cap \mho)$. This implies that
$$
\#\{1 \leq t \leq n; f^t(y) \in \hat P^{(k)}_j(x) \cap \mho \} \geq
$$
$$
\#\{1 \leq t \leq n; f^t(y) \in P^{(k)}_j(x) \}- 1.
$$
Hence
$$
\mu(\hat P^{(k)}_j(x) \cap \mho) \geq \mu(P^{(k)}_j(x) ) \geq \delta \mu(\hat P^{(k)}_j(x) ).
$$

\qed

Now define $ \La_\delta^0(i, l) $, for $\delta> 0$,  as the set of
points $x \in T^s$ such that for $k= k(i, l)$, we have
$$
\mu(P_j^{(k)}(x)) \geq \delta \mu(\hat P_j^{(k)}(x))
$$
  and
$$
\mu(P_j^{(k)}(x)) \geq  \delta \mu(\tilde P_j^{(k)}(x)),
$$
 for an infinite
sequence $\vsi(x)$ of values of  $j$.

Define $\La_\delta(i, l):= \La_\delta^0(i, l) \cap \mho(r_i, \rho_l)$.

The next lemma, a kind of Vitali's covering lemma, will be useful to estimate
the measure of $\mho^c$:

\begin{lemma} \label{levit}
Given a neighborhood $V$ of \  $\mho^c \cap
\La_\delta(i, l)$, there exist sequences $x_q \in \mho^c \cap \La_\delta(i, l)$, $(j_q), j_q \in \vsi(x_q)
\subset \natural$, $q= 1, 2, \dots$, such that
\begin{enumerate}
\item The sets $\hat P^{(k)}_{j_q}(x_q)$, $q \in \natural$ are disjoint and contained in $V$;

\item
$
\mu\big((\mho^c \cap \La_\delta(i, l)) \setminus \cup_{q \in \natural} \hat P^{(k)}_{j_q}(x_q)\big)= 0.
$
\end{enumerate}
\end{lemma}
\pf

By standard measure theoretical arguments, a translation $\tau: T^s
\to T^s$ can be found in such way that
$$
\mu(\tau(\cup \{\partial \hat A; A \in \SP^{(k)}_j, k \geq 1, j \geq
1\}))= 0;
$$
where $\partial \hat A$ is the boundary of $\hat A \in \hat
\SP^{(k)}_j$.

Denoting by $\SF$ the family of sets $P^{(k)}_j(x)$ with $x \in
\Lambda_\delta(i, l) \cap \mho^c$ and $j \in
\vsi(x)$. Take a sequence $A_u \in \SF$ satisfying:
\begin{enumerate}
\item $\hat A_u \subset V$, $\forall u \in \natural$, and $\mu(\hat
A_u \cap \hat A_e)= 0$, $\forall 1 \leq e< u$.

\item
\label{itemc}
$ \diam(A_u) = \max\{\diam(A); \hat A \subset V \text{ and }\mu(\hat
A \cap \hat A_e), \forall 1 \leq e< u\} $.
\end{enumerate}

Such properties imply that $\lim_{u \to + \infty} \diam(A_u)= 0$ and
\begin{equation}
\sum_{u} \mu(A_u)= \mu(\cup_u A_u) \leq 1. \label{eqsum}
\end{equation}

We claim that for $N \geq 1$
\begin{equation}
\big( \Lambda_\delta(i, l) \cap \mho^c \big) \setminus \cup_{u=
1}^N \ov {\hat A_u} \subset \cup_{u> N} \tilde A_u. \label{eqar}
\end{equation}

In fact, if $x \in \big( \Lambda_\delta(i, l) \cap \mho^c \big) \setminus \cup_{u= 1}^N \ov {\hat A_u}$, there exist $A \in \SF$ with $x \in A$ and
$$
\hat A \cap (\cup_{u= 1}^N \ov{\hat A_u})= \emptyset.
$$
Take $N_1> N$ such that $\hat A \cap \hat A_u= \emptyset$, $\forall 1 \leq u < N_1$
and $\hat A \cap \hat A_{N_1} \neq \emptyset$.
By item (\ref{itemc}) above, it follows that $\diam(\hat A) \leq \diam(\hat A_{N_1})$.
This implies that $\hat A \subset \tilde A_{N_1}$ and then
$$
x \in A \subset \tilde A_{N_1} \subset  \cup_{u> N} \tilde A_u,
$$
which concludes the proof of equation \ref{eqar}.
By such equation and our assumption that partition elements borders have zero measure,
we obtain that
$$
\mu\Big(
\big( \Lambda_\delta(i, l) \cap \mho^c \big) \setminus \cup_{u=
1}^N  {\hat A_u}
 \Big)=
\mu\Big(
\big( \Lambda_\delta(i, l) \cap \mho^c \big) \setminus \cup_{u=
1}^N \ov {\hat A_u}
 \Big) \leq
$$
$$
\mu\big(\cup_{u> N} \tilde A_u\big) \leq \sum_{u > N} \mu(\tilde A_u) \leq
\delta^{-1} \sum_{u > N} \mu(A_u).
$$
Due to eq. (\ref{eqsum}) the tail sum above goes to zero as $N \to +\infty$,
which implies
the lemma.

\qed

Lemmas  \ref{lerg} and \ref{levit} are the key ingredients in the

\begin{theorem} \label{teou}
Let $M$ be a compact Riemannian manifold and let $f:M \to M$ be a measurable Borelian map
such that $\SM_1(f) \neq \emptyset$.
Let $\mho \subset M$ and $\mho( r, \rho)$ be Borelian subsets of $M$,
where $r> 0$, $\rho> 1$.

Suppose that the points $x \in \mho(r, \rho)$ have the following
property: if $y \in B_{r'}(x)$ for some $0< r' \leq r$ and $f^t(y)
\in B_{r'}(x)$, for some $t> 0$ then there exist $0 \leq t_1 \leq
t$, such that $f^{t_1}(y) \in \ov{B_{\rho r'}(x)} \cap \mho$.
Suppose also that { $r_i> 0$, $\rho_l > 1$ are two monotone
sequences} converging respectively to 0 and $+\infty$, such that
$\cup_{i, l} \mho(r_i, \rho_l)= M$. Then, $\mho$ has total
probability with respect to the map $f$.

\end{theorem}
\pf
Consider  $\La_\delta^0(i, l)$ and $\La_\delta(i, l)= \La_\delta^0(i, l) \cap \mho(r_i, \rho_l)$
the same sets defined above in our text.

By Lemma \ref{leprp}, this implies that
$$
\mu(\La_\delta^0(i, l)) \geq 1- \delta (k^s + 3^s k^s).
$$
The last inequality implies that
$$
 \cup_{n= 1} ^{+\infty} \La_{1/n}(i, l)= \mho(r_i, \rho_l) \mod(0).
$$

Therefore, all we need is to prove that
$$
\mu(\mho^c \cap \La_\delta(i, l))= 0, \forall 0< \delta< 1,
$$
which implies that
$$
\mho \supset \mho( r_i, \rho_l) \mod(0) \Rightarrow \mho \supset \big(\cup_{i, l} \mho( r_i, \rho_l)\big) \mod(0) = M \mod(0),
$$
this last equality by hypothesis.
We will then have
$\mu(\mho)= \mu(M)= 1$, and the proof of Th. \ref{teou} will be completed.

Fix $(i, l)$ and $\delta> 0$. Let $V$ a neighborhood of  $\La_\delta(i, l)$.
 By lemmas
\ref{lerg} and \ref{levit} it follows that
$$
\mu(V) \geq \sum_q \mu\big(\hat  P^{(k)}_{j_q}(x_q)  \big) \geq
\frac{1}{1- \delta} \sum_q \mu\big(\hat  P^{(k)}_{j_q}(x_q) \cap \mho^c \big)=
$$
$$
\frac{1}{1- \delta}  \mu\Big(\big(\cup_q \hat  P^{(k)}_{j_q}(x_q) \big)\cap \mho^c \Big) \geq
\frac{1}{1- \delta}  \mu\big( \Lambda_\delta(i, l)\cap \mho^c \big).
$$

But if $\mu\big( \Lambda_\delta(i, l)\cap \mho^c \big)> 0$, one can take
$V$ satisfying
$$
\mu(V) < \frac{1}{1- \delta}  \mu\big( \Lambda_\delta(i, l)\cap \mho^c \big),
$$
contradicting the last inequality.
Hence $\mu\big( \Lambda_\delta(n, m)\cap \mho^c \big)= 0$.

\qed

Now, let us finish the proof of Proposition \ref{prop} (which implies Th. \ref{theop}).

\pf ({\bf Proposition \ref{prop}.})

Define  $\Sigma(\SU, \epsilon, r, \rho)$,
where $r> 0$, $\rho> 1$,
as the set of
points $x \in M$ such that  if $y \in B_{r'}(x)$ for some $0< r'
\leq r$ and $f^t(y) \in B_{r'}(t)$, for some $t> 0$ then there exist $0 \leq t_1< t_2
\leq t$, $g \in \SU$ and $z \in M$ such that $g= f$ on $M \setminus
B_\epsilon(f, x)$,
$$
g^{t_2- t_1}(z)= z, d(g^j(z), f^j(f^{t_1}(y)) \leq \epsilon, \forall
0 \leq j \leq t_2- t_1,
$$
and
$$
 f^{t_1}(y) \in \ov{B_{\rho r'}(x)}.
$$
(In particular, $f^{t_1}(y) \in \Sigma(\SU, \epsilon)$.)

It is easy to see that  $\Sigma(\SU, \epsilon, r, \rho)$ is a Borelian set.
Again, let $r_i> 0$ and $\rho_l> 1$ to be two monotone sequences converging respectively to $0$
 and $+\infty$.
We note that
Lemma \ref{lesomb} implies that
$$
M= \cup_{i \geq 1} \cup_{l \geq 1} \Sigma(\SU, \epsilon, r_i, \rho_l),
$$
for every neighborhood $\SU$ of $f$ and every $\epsilon> 0$.

So, taking $\mho= \Sigma(\SU, \epsilon)$ and $\mho(r_i, \rho_l)= \Sigma(\SU, \epsilon, r_i, \rho_l)$ in Th. \ref{teou}, we conclude that $\mu(\mho)= \mu(\Sigma(\SU, \epsilon))= 1$ for
all $f-$ergodic probability. By Ergodic Decomposition Theorem, this implies that $\Sigma(\SU, \epsilon)$ has total probability.

\qed

So far, we have proven the raw version of Ergodic Closing Lemma for
Endomorphisms (Th. \ref{theop}). The next lemma will be used in the
proof of the residual version of Ergodic Closing Lemma. We denote by
$\SM(M)$ the set of probabilities on $M$ endowed with the weak-*
topology.

\begin{lemma} \label{leres}
Let $f:M \to M$ be an endomorphism. Suppose that, for $x$ in a total probability set $S \subset M$,
given $\epsilon > 0$ and $\SU$ a neighborhood of $f$, there exists
$g_{x,  \epsilon} \in \SU$ and a $g_{x, \epsilon}-$periodic point
$p= p(x, \epsilon)$ which $\epsilon-$shadows $x$. Then, given any
ergodic measure $\mu \in \SM_1(f)$, there are $g_k \to f$ and
$g_k-$periodic points $p_k$ such that $\mu$
 is the limit of the sequence $(\mu_k)$ of $g_k-$ergodic measures respectivelly supported in
the orbit of $p_k$. Moreover, each $p_k$ can be taken to be a
hyperbolic periodic point for $g_k$.
\end{lemma}
\pf

Let us consider an $f$-ergodic probability $\mu$. We suppose,
without loss of generality, that $\mu$ is not supported in a
periodic orbit, otherwise there is nothing to prove. For a
$\mu-$typical point $x\in M$, we can assume that $x$ is recurrent
(by Poincar\'e's Recurrence Theorem), has the shadowing property as
in lemma's statement, and that
\begin{equation}
\frac 1 n \sum_{j= 0}^{n-1} \delta_{f^j(x)} \to_{_{\text{weak}-*}}
\mu, \label{eqfrac}
\end{equation}
as $n \to +\infty$ (this, by the Ergodic Decomposition Theorem).

Set $\epsilon_1= 1$ and $n_k> 0$ as the first return time of the
orbit of $x$ to $B(x, \epsilon_k)$, where $\epsilon_{k+ 1}:=
d(f^{n_{k}}(x), x)/2$, $\forall k \geq 1$.

Therefore, $n_k \to +\infty$ as $k \to +\infty$. By hypothesis, one
can take  a sequence of $g_k:= g_{x, \epsilon_k}$, with $g_k \to f$,
exhibiting  $g_k-$periodic points $(p_k)$ such that each $p_k$
$\epsilon_k/3$-shadows the orbit of $x$. In particular, the period
$t_{k+1}$ of $p_{k+1}$ is, at least, $n_k$ (otherwise, the orbit of
$x$ would return to $B(x, \epsilon_k)$ before $n_k$). So, $t_{k+1}
\geq n_k$ implies that $t_k \to +\infty$ as $k \to +\infty$, and (up
to take a subsequence) we can suppose that $t_k$ are distinct. Note
that slightly perturbing $g_k$ in the neighborhood of $p_k$, we can
suppose that $p_k$ is hyperbolic. Set $\mu_k$ as the $g_k$-ergodic
probability supported in the orbit of $p_k$. We will show that
$\mu_k \to_{_{\text{weak}-*}} \mu$ as $k \to +\infty$. From equation
\ref{eqfrac} we have that
$$
\nu_k= \frac 1 {t_k} \sum_{j= 0}^{t_k-1} \delta_{f^j(x)}
\to_{_{\text{weak}-*}} \mu,
$$
as $k \to +\infty$. Let $\alpha> 0$ and $\{\vr_1, \dots,\vr_s\}
\subset C^0(M)$ be given. All we need to see is that there exists
$k_0 \in \natural$ such that $\mu_k$ belongs to the neighborhood
$$
V_{\vr_1, \dots, \vr_s; \alpha} := \{ \nu \in \SM(M); |\int \vr_i
d\nu - \int \vr_i d\mu|< \alpha, \forall i= 1, \dots, s\},
$$
forall $k \geq k_0$.

 In fact, as $\vr_i$, $i= 1, \dots, s$ are uniformly continuous,
 take $\epsilon> 0$ such that $|\vr_i(y)- \vr_i(z)|<
\alpha/2$, $\forall i= 1, \dots, s$, $\forall y, z \in M$ such
that $d(y, z)< \epsilon$. Then, take $k_0$ such that $\epsilon_k<
\epsilon/2$,and $|\int \vr_i d\nu_k - \int \vr_i d\mu | < \alpha/2$,
$\forall k \geq k_0$, $\forall i= 1, \dots, s$. We conclude that
$$
|\int \vr_i d\mu_k - \int \vr_i d\mu| \leq |\int \vr_i d\mu_k - \int
\vr_i d\nu_k| + |\int \vr_i d\nu_k - \int \vr_i d\mu| <
$$
$$
\frac 1 {t_k} \sum_{j= 0}^{t_k-1} |\vr_i(f^j(x))- \vr_i(g_k^j(p_k))|
+ \alpha/2 \leq \alpha, \forall i= 1, \dots, s;
$$
which implies the lemma.

\qed

Now, we proceed with the proof of Th. \ref{theo1}, by deriving it
from Th. \ref{theop} and lemma \ref{leres} above. The arguments here
are basically the same as in Th. 4.2 in \cite{ABC}.

\pf ({\bf Th. \ref{theo1} - Ergodic Closing Lemma for Nonsingular Endomorphisms -
 Residual version.})

 For  $m \in \natural$ fixed, by standard transversality arguments,
the collection $\SK_m$ of  endomorphisms $f$ such that all periodic
points of $f$,  with period up to $m$ are hyperbolic is an open and
dense subset  $\NE(M)$.

So $\hat \SR:= \cap_{m= 1}^{+\infty} \SK_m$ is a residual set. Let
the set of probabilities $\SM(M)$ on $M$ to be endowed with the
weak-* topology and let $\kappa$ be the collection of compact
subsets of $\SM(M)$ endowed with Hausdorff distance. Given $f \in
\SR$, denote by $\SM_{per}(f)$ the set of $f-$ergodic measures
supported in $f-$periodic orbits.

Set $\Upsilon: \hat \SR \to \kappa$ given by
$$
\Upsilon(f)= \ov{\SM_{per}(f)}
$$
 Due to the
robustness of hyperbolic periodic points, such $\Upsilon$ is lower
semicontinuous. This implies that there is a residual subset $\SR
\subset \hat\SR$ whose elements are continuity points for
$\Upsilon$.

From now on, let $f \in \SR$. Let us prove that $\SM_1(f)$ is the
closed convex hull of $f-$ergodic measures supported in $f-$periodic
orbits. By Ergodic Decomposition Theorem, all we need to prove is
that any $f-$ergodic measure $\mu$ is in $\ov{\SM_{per}(f) }$.

By Lemma \ref{leres}, such measure $\mu$ is accumulated by $\mu_k
\in \SM_{per}(g_k)$, where $g_k \to f$ as $k \to +\infty$. As $\hat
\SR$ is residual, by means of a slight perturbation, we can suppose
that $g_k \in \hat \SR$ (as we construct $p_k$ to be hyperbolic in
the proof of that Lemma \ref{leres}, such $p_k$ persist under any
sufficiently small perturbation).

Since $f$ is a continuity point for $\Upsilon$, we have that
$\ov{\SM_{per}(g_k)} \to \ov{\SM_{per}(f)}$ as $k \to +\infty$, and
this implies that $\mu \in \ov{\SM_{per}(f)}$.

Therefore, $\ov{\SM_{per}(f)}$ contains all $f-$ergodic measures,
and by Ergodic Decomposition Theorem, we conclude that $\SM_1(f)$ is
the closed convex hull of $f$-ergodic measures supported in periodic
orbits.

\qed

\end{document}